\documentclass{cmslatex}
\usepackage[paperwidth=7in, paperheight=10in, margin=.875in]{geometry}
 \usepackage[backref,colorlinks,linkcolor=red,anchorcolor=green,citecolor=blue]{hyperref}
\usepackage{amsfonts,amssymb}
\usepackage{amsmath}
\usepackage{graphicx}
\usepackage{cite}
\usepackage{enumerate}
\usepackage{mathtools}
\renewcommand{\theequation}{\arabic{section}.\arabic{equation}}

\newcommand{\cM}{\mathcal{M}}

\newcommand{\cH}{\mathcal{H}}

\newcommand{\cK}{\mathcal{K}}
\newcommand{\cC}{\mathcal{C}}
\newcommand{\cJ}{\mathcal{J}}

\newcommand{\bN}{\mathsf{N}}

\newcommand{\yw}[1]{\textcolor{black}{#1}}

   \allowdisplaybreaks
\begin{document}
 \title{Second order ensemble Langevin method for sampling and inverse problems}

    
    \author{Ziming Liu\thanks{Department of Physics,
    Massachusetts Institute of Technology, Cambridge, Massachusetts 02139, zmliu@mit.edu.}
    \and Andrew Stuart\thanks{Applied and Computational Mathematics, California Institute of Technology, Pasadena, California 91125, astuart@caltech.edu} \and Yixuan Wang\thanks{Applied and Computational Mathematics, California Institute of Technology, Pasadena, California 91125, roywang@caltech.edu}}
    
    \pagestyle{myheadings} \markboth{Second order ensemble Langevin method}{Z Liu, A Stuart, Y Wang} 
    
    
    \date{\today}
    
    \maketitle
    
    \begin{abstract} We propose a sampling method based on an ensemble approximation of second order Langevin dynamics. {The log target density is appended with a quadratic term in an auxiliary momentum variable and 
damped-driven Hamiltonian dynamics, is introduced; the resulting stochastic
    differential equation is invariant to the Gibbs measure, with marginal on
    the position coordinates given by the target. A preconditioner based on 
covariance under the law of \yw{position coordinates under}
the dynamics does not change this invariance property, and is introduced 
to accelerate convergence to the Gibbs measure.}
        The resulting mean-field dynamics may be approximated by an ensemble 
method; this results in a gradient-free and affine-invariant stochastic 
dynamical system \yw{with desirable provably 
uniform convergence properties across
the class of all Gaussian targets.} Numerical results demonstrate the
potential of the method as basis for a numerical sampler in Bayesian inverse 
problems, \yw{beyond the Gaussian setting.}
    \end{abstract}
    
    \begin{keywords}  \yw{{Bayesian inverse problems; sampling; ensemble method; 
second order Langevin equation; Hybrid Monte Carlo; mean-field models; 
nonlinear Fokker-Planck equation}.}
    \end{keywords}

\begin{AMS} 62F15; 60H10; 65C30; 65C35
\end{AMS}
    \section{Introduction}
    \subsection{Set-Up}
    Consider sampling the density $$\pi(q)=\frac{1}{Z_q}\exp(-\Phi(q))\,,$$ where $\Phi: \mathbb{R}^{N} \to \mathbb{R}$ is termed the \emph{potential} function and $Z_q$  the \emph{normalization constant.}
    A broad family of problems can be cast into this formulation; the Bayesian approach to inverse problems provides a particular focus for our work \cite{kaipio2006statistical}. 
        The point of departure for the algorithms considered in this paper is the following mean-field
    Langevin equation:
\begin{equation}
\label{eqol}
\frac{\mathrm{d} q}{\mathrm{d} t}=-\cC(\rho) D \Phi(q)+\sqrt{2 \cC(\rho)} \frac{\mathrm{d} W}{\mathrm{d} t}\,,
\end{equation}
where $D$ \yw{denotes} the gradient operator, $W$ is an $N$-dimensional standard Brownian motion, $\rho$ is the density associated to the law of $q$, and $\cC(\rho)$ is the covariance under this density. This
constitutes a mean-field generalization \cite{garbuno2020interacting} of the standard Langevin
equation \cite{pavliotis2014stochastic}. Applying a particle approximation to the mean-field
model results in an interacting particle system, and coupled Langevin dynamics 
\cite{leimkuhler2018ensemble,garbuno2020interacting,nusken2019note,garbuno2020affine}. 
The benefit of preconditioning using the covariance is that it leads to
mixing rates independent of the problem, provably for quadratic $\Phi$ and
empirically beyond this setting \cite{garbuno2020interacting}, 
{because of the affine invariance \cite{goodman2010ensemble} of the
resulting algorithms \cite{garbuno2020affine}.}

    In order to {further} accelerate mixing and achieve sampling efficiency, 
    we introduce an auxiliary variable  $p\in \mathbb{R}^{N}$ and consider the Hamiltonian \begin{equation}
    \label{ham}
        \cH(z)=\frac{1}{2}\langle p,\cM^{-1}p\rangle+\Phi(q)\,,
    \end{equation} where the new state variable $z:=(q^\top, p^\top)^\top \in\mathbb{R}^{2N}$. Define a measure {via its density} on $\mathbb{R}^{2N}$ by
     
\begin{equation}
\label{eq:Hneed}
\Pi(z)=\frac{1}{Z_{q,p}}\exp(-\cH(z))\,,
\end{equation}
where  $Z_{q,p}$ is the normalization constant. The marginal distribution of $\Pi$ in the $q$ variable gives the desired distribution $\pi$, i.e. $\int \Pi(z) \mathrm{d}p=\pi(q)$.
    We now aim at sampling the joint distribution. To this end, consider the following underdamped Langevin dynamics in $\mathbb{R}^{2N}$:
\begin{equation}
\label{eqsol}
\frac{\mathrm{d} z}{\mathrm{d} t}=\cJ D {\cH}(z)-\mathcal{\cK} D {\cH}(z)+\sqrt{2 \mathcal{\cK}} \frac{\mathrm{d} W_0}{\mathrm{d} t}\,,
\end{equation}
with the choices 
\begin{equation}
\label{eqsol2}
\cJ=\begin{pmatrix}
0 & \cC \\
-\cC & 0
\end{pmatrix}\,,\quad
\mathcal{K}=\begin{pmatrix}
\mathcal{K}_{1} & 0 \\
0 & \mathcal{K}_{2}
\end{pmatrix}\,.
\end{equation}
Here $W_0$ is a  standard Brownian motion in $\mathbb{R}^{2N}$ with components $W', W \in \mathbb{R}^{N}$. Then we have the following, proved in Subsection
\ref{proofprop0}:

\vspace{0.1in}

\begin{proposition}
\label{prop.prop0}
{Assume that $\mathcal{K}_1$, $\mathcal{K}_2$ are symmetric and 
non-negative definite, and that $\mathcal{C}$ is symmetric positive definite. 
Assume further that $\mathcal{C}$, $\mathcal{K}$ and $\cM$ depend
on the law of $z$ under the dynamics defined by \eqref{eqsol} and \eqref{eqsol2},
but are are independent of $z$: all derivatives with respect to $z$
are zero. Then the Gibbs measure $\Pi(z)$ is invariant under the dynamics defined
by \eqref{eqsol}, \eqref{eqsol2}.}
\end{proposition}

\vspace{0.1in}

{In practice, to simulate from such a mean-field model, it will
be necessary to consider a particle approximation of the form
\begin{equation}
\label{eqsolp}
\frac{\mathrm{d} z^{(i)}}{\mathrm{d} t}={J(Z) D_{z} H(z^{(i)};Z)-K(Z) D_{z} H(z^{(i)};Z)}+\sqrt{2 K(Z)} \frac{\mathrm{d} W_0^{(i)}}{\mathrm{d} t}\,,
\end{equation}
for the set of $I$ particles $Z=\{ z^{(i)} \}_{i=1}^I,$ and where {$M(Z), K(Z), J(Z)$ are
appropriate empirical approximations of $\mathcal{M}(\rho), \mathcal{K}(\rho), \mathcal{J}(\rho)$ based on replacing $\rho$ by $\rho^I$ where}
$$\rho^I = {\frac{1}{I}}\sum_{i=1}^{I} \delta_{z^{(i)}}\,,$$}
and the Hamiltonian is given by
\begin{equation}
\label{eq:surelyneed}
H(z;Z)=\frac{1}{2}\langle p,M(Z)^{-1}p\rangle+\Phi(q)\,.
\end{equation} 
Thus
$$D_z H(z;Z)=\bigl(D\Phi(q)^\top, (M(Z)^{-1}p)^\top\bigr)^\top.$$
\yw{Note that $H$ is the appropriate finite particle approximation of
$\mathcal{H},$ given the particle approximation $M$ of $\mathcal{M}.$ The 
dependence of $\mathcal{H}$ on the law of $z$ has been replaced by
dependence on the collection of particles $Z$.\footnote{The reader is asked
to note that collection of particles $Z$ is different from normalization
constant $Z_{q,p}$ appearing in \eqref{eq:Hneed}.}}

\vspace{0.1in}

\begin{remark}
\label{rem:12}
\yw{Unlike \eqref{eqsol}, the equation \eqref{eqsolp} is no longer 
a damped-driven Hamiltonian system; this is because of the dependence
of the Hamiltonian on the particle positions $Z$, through the mass
matrix. Furthermore, its marginal on any 
coordinate $q^{(i)}$ does not necessarily
preserve the desired target measure under the dynamics. However we
expect it to do so \emph{approximately} when $I$ is large. This
justifies the use of algorithms based on \eqref{eqsolp}.}
\end{remark}

\vspace{0.1in}

{In this paper we will concentrate on a specific choice of mean-field operators within
the above general construction, which we now describe.} Let $\cC_q(\rho)$ denote the $q$-marginal in the covariance 
under the law of \eqref{eqsol}. We make the choices
{$\mathcal{K}_{1}=0$, $\cC=\cM=\cC_q(\rho)$, and $\cK_2=\gamma\cC_q(\rho)$, for a scalar damping parameter $\gamma>0$.} Then the underdamped Langevin dynamics yields
\begin{equation}
\label{eqmom}
\begin{aligned}
\frac{\mathrm{d} q}{\mathrm{d} t}&= p\,, \\
\frac{\mathrm{d} p}{\mathrm{d} t}&=-\cC_q(\rho)D \Phi(q)- {\gamma p+\sqrt{2\gamma \cC_q(\rho)} \frac{\mathrm{d} W}{\mathrm{d} t}}\,.
\end{aligned}
\end{equation}

{To implement a particle approximation of the}
mean-field dynamics \eqref{eqmom} {we introduce  particles in the form} $z^{(i)}(t)=\bigl( (q^{(i)}(t))^\top, (p^{(i)}(t))^\top\bigr)$  and
we use the {ensemble covariance and mean approximations}
\begin{subequations}
\label{eq:ecov}
\begin{align}
    \cC_q(\rho)\approx C_q(Z)&\coloneqq\frac{1}{I} \sum_{i=1}^{I}\left(q^{(i)}-\bar{q}\right) \otimes\left(q^{(i)}-\bar{q}\right)\,,\\
    \bar{q}&:=\frac{1}{I} \sum_{i=1}^{I}q^{(i)}\,.
\end{align}
\end{subequations}
 In order to 
{obtain} affine invariance, we take the generalized square root of the ensemble covariance $C_q(Z)$, {similarly} to \cite{garbuno2020affine}.  We introduce the $N\times I$ matrix $$Q \coloneqq\left(q^{(1)}-\bar{q},q^{(2)}-\bar{q},\cdots,q^{(I)}-\bar{q}\right)\,,$$ which allows us to  
define the empirical covariance and generalized (nonsymmetric) square root  via
$${C}_q(Z)=\frac{1}{{I}} QQ^T\,,\quad
\sqrt{{C}_q(Z)}\coloneqq\frac{1}{\sqrt{I}} Q\,.$$ 
Now with $I$ independent standard Brownian motions $\{ W^{(i)} \}_{i=1}^I\in \mathbb{R}^{I}$, a natural
particle approximation of \eqref{eqmom} is, for $i=1, \dots, I$,
\begin{equation}
\label{ensol}
\begin{aligned}
\frac{\mathrm{d} q^{(i)}}{\mathrm{d} t}&= p^{(i)}\,, \\
\frac{\mathrm{d} p^{(i)}}{\mathrm{d} t}&=-C_q(Z)D \Phi(q^{(i)})- {\gamma p^{(i)}+\sqrt{2 \gamma C_q(Z)} \frac{\mathrm{d} W^{(i)}}{\mathrm{d} t}}\,.
\end{aligned}
\end{equation}

\yw{In subsequent sections we will 
employ an ensemble approximation of $C_q(Z)D \Phi(\cdot)$, as in \cite{garbuno2020interacting}, thereby
avoiding the need to compute adjoints of the forward model; we note also that in the linear case this approximation is exact. We will show that the resulting interacting
particle system has the potential to provide accurate derivative-free inference
for {certain classes of} inverse problems.}

\yw{In the remainder of this section we provide a literature review, we
highlight our contributions, and we outline the structure of the paper.}


\subsection{Literature Review}
\label{ssec:lit}
The overdamped Langevin equation is the canonical SDE that is invariant with respect to \yw{a given target density}. Many sampling algorithms are built upon this idea, and in particular, it is shown to govern a large class of Monte Carlo Markov Chain (MCMC) methods; see \cite{roberts2001optimal,ottobre2011asymptotic,ma2015complete}.
To enhance mixing and accelerate convergence, a second-order method, Hybrid Monte Carlo (HMC, also referred to as Hamiltonian Monte Carlo) \cite{duane1987hybrid,neal1994improved} has been proposed, leading
to underdamped Langevin dynamics. There have been many attempts to justify the empirically
observed fast convergence speed of second-order methods in comparison with first-order methods \cite{betancourt2017conceptual}. Recently a quantitative convergence rate is established in \cite{cao2019explicit}, showing that the underdamped Langevin dynamics converges faster than the overdamped Langevin dynamics when the log of the target $\pi$ has a small Poincar\'e constant; see also \cite{eberle2019couplings}. 

\yw{The idea of introducing preconditioners within the context of interacting particle systems used for sampling is developed in  \cite{leimkuhler2018ensemble}.} Preconditioning via ensemble covariance {is} shown to boost convergence and numerical performance \cite{garbuno2020interacting}. Other choices of preconditioners in sampling will lead to different forms of SDEs and associated Fokker-Planck equations with different structures, which can result in different sampling methods effective in different scenarios; see for example \cite{lindsey2021ensemble}.
Affine invariance is introduced in \cite{goodman2010ensemble} where it is argued that
this property leads to desirable convergence properties for interacting
particle systems used for sampling: methods that satisfy affine invariance are invariant under an affine change of coordinates, and are thus uniformly effective for problems that can be rendered  well-conditioned under an affine 
transformation; \yw{in particular for the sampling of the class of
all Gaussian measures.} An affine invariant version of the mean-field 
underdamped Langevin dynamics of \cite{garbuno2020interacting} 
is proposed in \cite{nusken2019note,garbuno2020affine}.

     {Kalman methods have shown wide success in state estimation problems since their
     introduction by Evensen; see \cite{evensen2009data} for an overview of the field,
     and the papers \cite{reich2011dynamical,iglesias2013ensemble} for discussion of
     their use in inverse problems. Using the ensemble covariance as a preconditioner  leads to affine invariant \cite{garbuno2020affine} and gradient-free 
     \cite{garbuno2020interacting} approximations of Langevin dynamics; this is desirable in practical computations in view of the intractability of derivatives in many large-scale models arising in science and engineering \cite{cleary2021calibrate,haber2018never,kovachki2019ensemble,huang2021unscented}. 
     See \cite{bergemann2010mollified,schillings2017analysis} for analysis of these methods and, in the context of continuous data-assimilation, see \cite{del2017stability,bergemann2012ensemble,taghvaei2018kalman}. {There are other
     derivative-free methods that can be derived from the mean-field perspective,
     and in particular consensus-based methods show promise for optimization
     \cite{carrillo2018analytical} and have recently been extended to sampling in
     \cite{carrillo2022consensus}.}}

     Recent work has established the convergence of ensemble preconditioning methods to mean field limits; see for example \cite{ding2021ensemble}. For other works on rigorous derivation of mean field limits of interacting particle systems, see \cite{sznitman1991topics,carrillo2010particle,jabin2017mean}. For underpinning theory of Hamiltonian-based sampling, see
     \cite{bou2017randomized,bou2018geometric,livingstone2019geometric,betancourt2017geometric}.

\subsection{Our Contributions}
\label{ssec:con}
The following contributions are made in this paper:
\begin{enumerate}
    \item We introduce an underdamped second order mean field Langevin dynamics, with a covariance-based preconditioner. 
    \item {In the case of {Bayesian inverse problems defined by}
    a linear forward map, we show that 
    that {this mean field model} preserves Gaussian distributions under time evolution and, if initialized at
    a Gaussian, converges to the desired target at {a} rate independent of the linear map.}
    \item \yw{We introduce finite particle approximations of the mean field model, resulting in an affine invariant method.}
    \item For {Bayesian inverse problems}, we introduce a gradient-free approximation of the algorithm, based on ensemble Kalman methodology.
    \item {In the context of {Bayesian inverse problems} we provide numerical examples to demonstrate that the algorithm resulting
    from the previous considerations has desirable sampling properties.}
\end{enumerate}

In Section \ref{2}, we introduce the inverse problems context that motivates us.
In Section \ref{3} we discuss the equilibrium distribution of the mean field
model.
Section \ref{EKA} introduces the ensemble Kalman approximation of the
finite particle system; {and in that section we also
demonstrate affine invariance of the resulting method.} Section \ref{2l} presents analysis of the finite particle system in the case of linear inverse problems, where the ensemble Kalman approximation is exact; we demonstrate
that the relaxation time to equilibrium is independent of the specific linear inverse
problem considered, a consequence of affine invariance. In Section \ref{4} we provide numerical results which demonstrate the efficiency and potential value of our method, 
and in Section \ref{5} we draw conclusions. Proofs of the propositions are given in the appendix, Section \ref{6}.

\section{Inverse Problem}\label{2}

Consider the Bayesian inverse problem of finding $q$ from an observation $y$ determined by the forward model 
$$y=\mathcal{G}(q)+\eta\,.$$
{Here $\mathcal{G}:\mathbb{R}^N \to \mathbb{R}^J$ is a {(in general)} nonlinear forward map. We assume a prior zero-mean Gaussian $\pi_{0}=\mathsf{N}(0, \Gamma_0)$ on unknown $q$ and assume that
\yw{the random variable} $\eta \sim \mathsf{N}(0, \Gamma)$, representing measurement error, is independent of the prior on $q$. 
We also assume that $\Gamma, \Gamma_0$ are positive definite.} Then by Bayes rule, the posterior density that we aim to sample is given by\footnote{{In what follows $\|\cdot\|_{C}=\|C^{-\frac12}\cdot\|$, with analogous notation for the
inducing inner-product, for any positive definite covariance $C$ and for $\|\cdot\|$ the Euclidean norm.}}
\yw{$$\pi(q) \propto \exp \Bigl(-\frac{1}{2}\|y-\mathcal{G}(q)\|_{\Gamma}^{2}\Bigr) \pi_{0}(q) \propto {\rm exp}(-\Phi(q))\,,$$
where potential function $\Phi(q)$ has the following form:}
\begin{equation}
\label{invf}
\Phi(q)=\frac{1}{2}\|y-\mathcal{G}(q)\|_{\Gamma}^{2}+\frac{1}{2}\|q\|_{\Gamma_0}^{2}\,.
\end{equation}
 
{In the linear case when $\mathcal{G}(q)=Aq$, $\Phi(q)$ is quadratic
 and the gradient $D\Phi(q)$ can be written as a linear function:
\begin{subequations}
\label{eq:quad}
\begin{align}
\Phi(q)&=\frac{1}{2}\|y-Aq\|_{\Gamma}^{2}+\frac{1}{2}\|q\|_{\Gamma_0}^{2}\,,\\
D\Phi(q)&=B^{-1}q-c\,,\\
B&=(A^T\Gamma^{-1} A+\Gamma_0^{-1})^{-1}\,,\quad
c=A^T\Gamma^{-1} y\,.
\end{align}
\end{subequations}
{In this linear setting,} the posterior distribution $\pi(q)$ is the Gaussian $\bN(Bc,B)$.}

\section{\yw{Equilibrium Distributions for the Mean Field Fokker-Planck Equation}}
\label{3}
{The mean-field underdamped Langevin equation \eqref{eqsol} has an associated 
nonlinear and nonlocal Fokker-Planck equation giving the evolution of the law of particles
$z(t)$, denoted $\rho(z,t)$. The equation for this law is
{(see proof in subsection \ref{proofprop0}.)}
\begin{equation}
\label{fp}
\partial_t \rho=\nabla \cdot \bigl((\mathcal{K}-J)(\rho\nabla {\cH}+\nabla\rho)\bigr)\,.
\end{equation}
{By Proposition \ref{prop.prop0} this Fokker-Planck equation has $\Pi(z)$ as its equilibrium; this follows as for standard linear Fokker-Planck equations \cite{pavliotis2014stochastic} since the dependence of the sample paths on $\rho$ involves only the mean and covariance; see the proof in subsection \ref{proofprop0}
and see also \cite{duncan2017using, graham1977covariant,jiang2004mathematical,risken1996fokker,pavliotis2014stochastic} for discussions on how to derive Fokker-Planck equations with a prescribed stationary distribution. 
In the specific case \eqref{eqmom}, recall that $\cJ$ and $\mathcal{K}$ are 
given by}}
\begin{equation} 
\cJ=\begin{pmatrix}
0 & \cC_q(\rho) \\
-\cC_q(\rho) & 0
\end{pmatrix}\,,\quad
\mathcal{K}=\begin{pmatrix}
0 & 0 \\
0 & {\gamma}\cC_q(\rho)
\end{pmatrix}\,.
\end{equation}
\yw{These choices satisfy the assumption of Proposition \ref{prop.prop0}.
We approximate \eqref{eqmom} by the interacting particle system \eqref{ensol}.
In this context, we note Remark \ref{rem:12} to motivate computational methods
based on integrating \eqref{ensol}.}


\section{Ensemble Kalman Approximation}\label{EKA}
  
\subsection{Derivatives Via Differences}
\label{dvd}

{We now make the ensemble Kalman approximation to approximate the gradient term by differences, as in \cite{garbuno2020interacting}:
\[D \mathcal{G}\left(q^{(i)}\right)\left(q^{(k)}-\bar{q}\right) \approx\left(\mathcal{G}\left(q^{(k)}\right)-\bar{\mathcal{G}}\right)\,,\]
where $\bar{\mathcal{G}}:=\frac{1}{I} \sum_{k=1}^{I} \mathcal{G}\left(q^{(k)}\right)$. 
Invoking this approximation within \eqref{ensol}, using the specific form
\eqref{invf} of $\Phi$,
yields the following system of interacting particles in $\mathbb{R}^N$, for $i=1, \dots, I:$}
\begin{equation}
    \label{encgf}
\begin{aligned}
\dot{q}^{(i)}&= p^{(i)}\,,  \\
\dot{p}^{(i)}&=-C_q(Z)\Gamma_0^{-1} {q}^{(i)}
-\frac{1}{I} \sum_{k=1}^{I}\langle  \mathcal{G}(q^{(k)})-\bar{\mathcal{G}}, \mathcal{G}(q^{(i)})-y\rangle_{\Gamma} q^{(k)}- {\gamma}p^{(i)}+\sqrt{2{\gamma} C_q(Z)} \dot{W}^{(i)}\,.
\end{aligned}
\end{equation}  
{We will use this system as the basis of all our numerical experiments.}

\subsection{Affine Invariance}
\label{2a}

In this subsection, we show the affine invariance property \cite{goodman2010ensemble,garbuno2020affine,leimkuhler2018ensemble}  for the Fokker-Planck equations in the mean-field regime \eqref{fp}, for the particle equation in the mean-field regime \eqref{eqmom}, for the ensemble approximation \eqref{ensol}, and the gradient-free approximation  \eqref{encgf}. For simplicity of presentation, we only state the results in the case of ensemble approximation, and the mean-field case is a straightforward analogy upon dropping all of the particle superscripts.

\vspace{0.1in}

\begin{definition}[Affine invariance for particle formulation]
\label{af}
We say a particle formulation is affine invariant, if under all affine 
transformations of the form 
\begin{equation}
\label{eq:aft}
q^{(i)}=A v^{(i)}+b\,,\quad p^{(i)}=A u^{(i)}\,,
\end{equation}
the equations on the transformed particle systems are given by the same equations with $q^{(i)}$, $p^{(i)}$ replaced by $v^{(i)}$, $u^{(i)}$ respectively, and with potential $\Phi$ replaced by $\tilde{\Phi}$ via \[\tilde{\Phi}(v^{(i)})=\Phi(q^{(i)})=\Phi(Av^{(i)}+b)\,.\]
Here $A$ is any invertible matrix and $b$ is a vector.
\end{definition}

\vspace{0.1in}

\begin{definition}[Affine invariance for Fokker-Planck equation]
\label{af1}
We say a Fokker-Planck equation is affine invariant, if under all affine transformations of the form \[q^{(i)}=A v^{(i)}+b\,,\quad p^{(i)}=A u^{(i)}\,,\]
the equations on the pushforward PDF $\tilde{\rho}^I$ are given by the same equation on ${\rho}^I$ with $q^{(i)}$, $p^{(i)}$ replaced by $v^{(i)}$, $u^{(i)}$ respectively, and with Hamiltonian $H$ replaced by $\tilde{H}$ via
\[ \tilde{H}(v^{(i)}, u^{(i)})=H(q^{(i)}, p^{(i)})=H(Av^{(i)}+b,Au^{(i)})\,.\]
Here $A$ is any invertible matrix and $b$ is a vector.
\end{definition}

\vspace{0.1in}

The key dynamical systems introduced in this paper are affine invariant:

\vspace{0.1in}

\begin{proposition}
\label{aft}
The particle formulations \eqref{eqmom}, \eqref{ensol} and \eqref{encgf} are affine invariant. The Fokker-Planck equation \eqref{fp} is also affine invariant.
\end{proposition}

\vspace{0.1in}

{We defer the proof to Subsection \ref{proofaf}. The significance of affine invariance
is that it implies that the rate of convergence is preserved under affine transformations.
The proposed methodology is thus uniformly effective for problems that become well-conditioned under an affine transformation.
Proposition \ref{afli}, which follows in the next section, illustrates this property
in the setting of linear forward map $\mathcal{G}(\cdot).$}

\vspace{0.1in}

\begin{remark}
The affine invariance of the methodology introduced
in \cite{leimkuhler2018ensemble} involves a definition different
from that in Definition \ref{af}. In particular \eqref{eq:aft}
is replaced by 
\begin{equation}
\label{eq:aft2}
q^{(i)}=A v^{(i)}+b\,,\quad p^{(i)}=u^{(i)}\,,
\end{equation} 
\end{remark}
\section{Mean Field Model For Linear Inverse Problems}
\label{2l}

{We consider the mean field SDE \eqref{eqmom}
in the linear inverse problem setting of Section \ref{2} with $\mathcal{G}(q)=Aq;$ thus
\eqref{eq:quad} holds.  We note that $B$ in \eqref{eq:quad}
is both well-defined and symmetric positive definite since
$\Gamma_0, \Gamma$ are assumed to be symmetric positive definite.
\yw{The two Propositions \ref{afli}, \ref{prop3}
demonstrate problem-independent rates of convergence, across the
set of all linear Gaussian inverse problems; this is a consequence} 
of affine invariance which in turn is a consequence
of our choice of preconditioned mean field system.

In the setting of the linear inverse problem, the mean field model \eqref{eqmom} reduces to
\begin{equation}
\label{linearp}
\begin{aligned}
\dot{q}&= p\,, \\
\dot{p}&=-\cC_q(\rho)(B^{-1}q-c) - {\gamma}{p}+\sqrt{2 {\gamma}\cC_q(\rho)}\dot{W}\,.
\end{aligned}
\end{equation}
We prove the following result about this system in Subsection \ref{proof3.3}:}

\vspace{0.1in}

\begin{proposition} 
\label{afli}
Write the mean $m(\rho)$ and the covariance $\cC(\rho)$ of the law $\rho(z)$  of particles in equation \eqref{linearp} in the block form
$$m(\rho)=\begin{pmatrix}
m_q(\rho) \\
m_p(\rho)
\end{pmatrix}\,, \quad \cC(\rho)=\begin{pmatrix}
\cC_q(\rho) & \cC_{q,p}(\rho) \\
\cC_{q,p}^T(\rho) & \cC_p(\rho)
\end{pmatrix}\,.$$ The evolution of the mean and covariance is prescribed by the following system of ODEs:
\begin{equation}
\label{vare}
\begin{aligned}
\dot{m_q}&= m_p\,, \\
\dot{m_p}&=-\cC_q(B^{-1}m_q-c) - {\gamma}m_p\,,\\
\dot{\cC_q}&=\cC_{q,p}+\cC_{q,p}^T\,,\\
\dot{\cC_{p}}&=-\cC_{q} B^{-1}\cC_{q,p}-(\cC_{q} B^{-1}\cC_{q,p})^T- 2{\gamma}\cC_{p}+2{\gamma}\cC_{q}\,,\\
\dot{\cC}_{q,p}&=- {\gamma}\cC_{q,p} -\cC_{q} B^{-1}\cC_{q} +\cC_{p}\,.
\end{aligned}
\end{equation}
{The unique steady solution with positive definite covariance is the Gibbs
measure $m_q=Bc, m_p=0, C_{q,p}=0, C_q=C_p=B;$ the marginal on $q$ gives the
solution of the linear Gaussian Bayesian inverse problem.
All other steady state solutions have degenerate covariance, are unstable, and take the
form $m_q=B(c+m), m_p=0, C_{q,p}=0, C_q=C_p=B^{1/2}XB^{1/2}$ for a projection matrix $X$ 
 and $m$ in the nullspace of $C_q$.
 The equilibrium point with positive definite covariance is hyperbolic
and linearly stable and hence its basin of attraction is an open set, containing
the equilibrium point itself, in the set of all mean vectors and 
positive definite
covariances. Furthermore, the mean and covariance converge to this equilibrium,
from all points in its basin of attraction, with a speed independent of $B$ and $c$.}  
\end{proposition} 
\vspace{0.1in}
\begin{remark}
\label{r:5.2}
{Proposition \ref{afli} demonstrates that the convergence speed to the
non-degenerate equilibrium point is independent of the specific linear
inverse problem to be solved; the rate does, however, depend on $\gamma$.
Analysis of the linear stability problem shows that $\gamma\approx1.83$ gives the best local convergence rate; see Remark \ref{optg}. In the case where mean field 
preconditioning is not used, the optimal choice of $\gamma$ for underdamped
Langevin dynamics depends on the linear inverse problem being solved, 
and can be identified explicitly in the scalar setting \cite{pavliotis2014stochastic}; for analysis in the non-Gaussian setting see \cite{chak2021optimal}. Motivated by analogies with the work in
\cite{pavliotis2014stochastic,chak2021optimal} we expect the optimal choice
of $\gamma$ to be problem dependent in the nonlinear case, but motivated
by our analysis in the linear setting we expect to see a good choice
which is not too small or too large and can be identified by straightforward
optimization.}
\end{remark}
\vspace{0.1in}

Now we show that for Gaussian initial data, the solution remains Gaussian and is thus determined by the evolution of the mean and covariance via the ODE \eqref{vare}.
We prove this by investigating the mean field Fokker-Planck equation in the linear case. 
The evolution in time of the law $\rho(z)$ of equation \eqref{linearp} is governed
by the equation
\begin{equation}
\label{linearfp}
\partial_t \rho=\nabla\cdot\left(
\begin{pmatrix}
-p \\
\cC_q(\rho)(B^{-1}q-c) + {\gamma}p
\end{pmatrix}\rho+\begin{pmatrix}
0 & -\cC_q(\rho) \\
\cC_q(\rho) & {\gamma}\cC_q(\rho)
\end{pmatrix}\nabla\rho\right)\,.
\end{equation}

\yw{We prove the following result in Subsection \ref{proof3}; by
virtue of Proposition \ref{afli} it establishes
that the eigenvalues that determine the
local stability of the posterior Gaussian density are
the same across all linear Gaussian inverse problems:} 

\vspace{0.1in}

\begin{proposition}
\label{prop3}
{Let $m(t)$, $C(t)$ solve the ODE \eqref{vare} with initial conditions $m_0$ and $C_0$.
Assume that $\rho_0$ is a Gaussian distribution, so that
\[\rho_{0}(z):=\frac{1}{(2 \pi)^{N}}\left(\operatorname{det} {C}_{0}\right)^{-1 / 2} \exp \left(-\frac{1}{2}\left\|z-m_{0}\right\|_{{C}_{0}}^{2}\right)\] 
with mean $m_0$ and covariance $C_0$. Then the Gaussian profile \[\rho(t,z):=\frac{1}{(2 \pi)^{N}}\left(\operatorname{det} {C}({t})\right)^{-1 / 2} \exp \left(-\frac{1}{2}\left\|z-m({t})\right\|_{{C}({t})}^{2}\right)\] solves the Fokker-Planck equation \eqref{linearfp} with initial condition $\rho(0,z)=\rho_{0}(z)$.}
\end{proposition}

\vspace{0.1in}

\section{Numerical Results}
\label{4}

{We introduce, and study, a numerical method for sampling the Bayesian inverse problem of Section \ref{2}; the method is based on numerical time  stepping of the interacting particle system \eqref{encgf}.}
{In this section, we demonstrate that the proposed sampler, which we refer to as
\emph{ensemble Kalman hybrid Monte Carlo} (EKHMC) can effectively approximate posterior distributions for two widely studied inverse problem test cases. We compare EKHMC with its first-order version EKS \cite{garbuno2020interacting} and a gold standard MCMC \cite{brooks2011handbook}. EKHMC inherits two major advantages of EKS: (1) exact gradients are not required (i.e., derivative-free); (2) the ensemble can faithfully approximate the spread of the posterior distribution, rather than collapse to a single point as happens with the
basic EKI algorithm \cite{schillings2017analysis}.
Furthermore, we show empirically that EKHMC can obtain samples of  similar quality to EKS, and has faster convergence than EKS. We detail our numerical time-stepping
scheme in the first subsection,
before studying two examples (one low dimensional, one a PDE inverse problem for a field) in
the subsequent subsections.}

\subsection{Time Integration Schemes}

We employ a splitting method to integrate the stochastic {dynamical} system
given by equation~(\ref{encgf}): \yw{the first capturing the finite
particle approximation of the Hamiltonian evolution,} 
and the second capturing an OU process in momentum space.  
\yw{The Hamiltonian evolution follows the equation}:
\begin{equation}\label{exp:Ham}
\begin{aligned}
\dot{q}^{(i)}&= p^{(i)}\,,  \\
\dot{p}^{(i)}&= F_{H}:= -C_q(Z)\Gamma_0^{-1} {q}^{(i)}
-\frac{1}{I} \sum_{k=1}^{I}\langle  \mathcal{G}(q^{(k)})-\bar{\mathcal{G}}, \mathcal{G}(q^{(i)})-y\rangle_{\Gamma} q^{(k)}\,.
\end{aligned}
\end{equation}
The OU process follows the equation:
\begin{equation}\label{exp:OU}
\begin{aligned}
\dot{q}^{(i)}&= 0\,,  \\
\dot{p}^{(i)}&=-\gamma p^{(i)} +\sqrt{2\gamma C_q(Z)} \dot{W}^{(i)}\,.
\end{aligned}
\end{equation}  
\yw{We implement a symplectic Euler  integrator \cite{sanz2018numerical} for 
the part of the particle system arising from approximation of
the Hamiltonian contribution
to the damped-driven mean-field equations. That is, we take a half 
step $\epsilon/2$ of momentum updates, then a full step $\epsilon$ of position updates, and finally a half step $\epsilon/2$ of momentum updates.
With the ensemble approximation the system is only approximately Hamiltonian; 
the splitting used in symplectic Euler is still well-defined, however. And we
also expect it to perform well in the large particle limit because
the mean-field limit is itself Hamiltonian.} {Let $Z_j$ be the collection of
all position and momentum particles at time $j$: $\{q_{j}^{(i)},p_{j}^{(i)}\}_{i=1}^I.$ Starting from time $j$ this
symplectic Euler integration gives map $Z_{j} \mapsto \hat{Z}_{j}.$
We set $q_{j+1}^{(i)}$ to be the $i^{\rm th}$ position coordinates of $\hat{Z}_{j}$
and then update the momentum coordinates using the OU process \yw{which provides
the damped-driven component of the mean-field limiting process.} 
The damping coefficient $\gamma>0$ is treated as a hyperparameter of EKHMC. 
{Vector $Z$ is set at the value given 
by output of the preceding symplectic Euler integrator, denoted by $\hat{Z}_{j}$.}
The update of the $i^{\rm th}$ momentum coordinate, \yw{given by 
solving the OU process exactly in law,} is then
\begin{equation}
\begin{aligned}
    & \tilde{p}^{(i)}_j={\rm exp}(-\gamma\epsilon)\hat{p}^{(i)}_j\,, \\
    & p^{(i)}_{j+1} = \tilde{p}^{(i)}_j+\eta,\quad \eta\sim \mathsf{N}
    \Bigl(0,\bigl(1-{\rm exp}(-2\gamma\epsilon)\bigr)C_q(\hat{Z}_{j})\Bigr)\, ,
\end{aligned}
\end{equation}
{where $\hat{p}^{(i)}_j$ are the momentum coordinates from $\hat{Z}_{j}$.}
Within the OU process the damping coefficient $\gamma>0$ is treated as a hyperparameter of EKHMC.

\yw{It is posible to consider use of a Metropolis-Hastings (Metropolization)
step to correct the dynamics; however because the underlying continuous
time system is not (for finite number of particles) invariant
with respect to the target, doing so would be very complicated
and so we do not pursue this. Aside from invariance,
Metropolization also imparts stability of the integrator and we may address
this in different ways. Indeed, similarly} 
to~\cite{garbuno2020interacting}, and as there 
with the goal of improving stability
and convergence speed towards posterior distribution, we implement an adaptive step size, i.e., the true step size is rescaled by the magnitude of the ``force field" $F_H$ (defined in equation~(\ref{exp:Ham})):
\begin{equation}
    \tilde{\epsilon}=\frac{\epsilon}{a|F_H|+1}\,.
\end{equation}

\subsection{Low Dimensional Parameter Space}

We follow the example presented in Section 4.3 of the paper~\cite{garbuno2020interacting}. We start by defining the forward map which is given by the one-dimensional elliptic boundary value problem
\begin{subequations}
\begin{align}
    &-\frac{d}{dx}\Bigl({\rm exp}(u_1)\frac{d}{dx}p(x)\Bigr)=1\,,\quad x\in(0,1)\,,\\
&\quad\quad\quad p(0)=0\,, \quad p(1)=u_2\,.
\end{align}
\end{subequations}
The solution is given explicitly by
\begin{equation}
    p(x) = u_2x+{\rm exp}(-u_1)\Bigl(-\frac{x^2}{2}+\frac{x}{2}\Bigr)\,,
\end{equation}
The forward model operator $\mathcal{G}$ is then defined by
\begin{equation}
    \mathcal{G}(u) = 
    \begin{pmatrix}
    p(x_1)\\
    p(x_2)
    \end{pmatrix}\,.
\end{equation}
Here $u=\bigl(u_1,u_2\bigr)^T$ is a constant vector that we want to find and we assume that we are given noise measurements $y$ of $p(\cdot)$ at locations $x_1=0.25$ and $x_2=0.75$. Parameters are chosen according to~\cite{garbuno2020interacting}, but we summarize them here for completeness:

\begin{itemize}
    \item noise $\eta\sim\mathsf{N}(0,\Gamma)$, $\Gamma=0.1^2I_2$;
    \item prior $\pi_0=\mathsf{N}(0,\Gamma_0)$, $\Gamma_0=\sigma^2I_2$, $\sigma=10;$
    \item measurement $y=(27.5,79.7)^T;$
    \item number of particles $I=10^3;$
    \item initialization: $(u_1,u_2)\sim \mathsf{N}(-3.5,0.1^2)\times \mathsf{U}(70,110).$
\end{itemize}
Here $\mathsf{U}$ is the uniform distribution.

\begin{figure}
    \centering
    \includegraphics[width=1.0\linewidth]{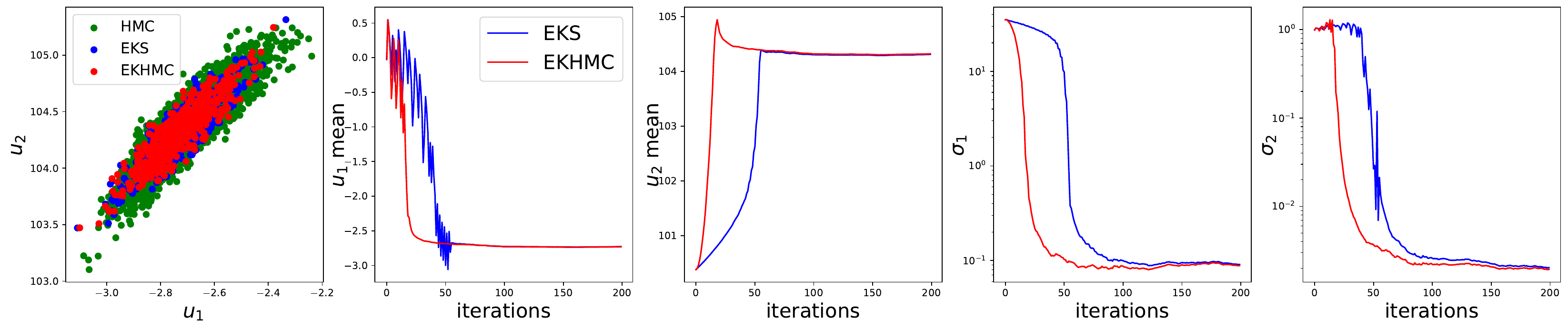}
    \caption{\yw{The low dimensional parameter space example. From left to right: samples; mean $u_1$; mean $u_2$; the first singular value $\sigma_1$; the second singular value $\sigma_2$.}}
    \label{fig:2d}
\end{figure}

\begin{figure}
    \centering
    \includegraphics[width=1\linewidth]{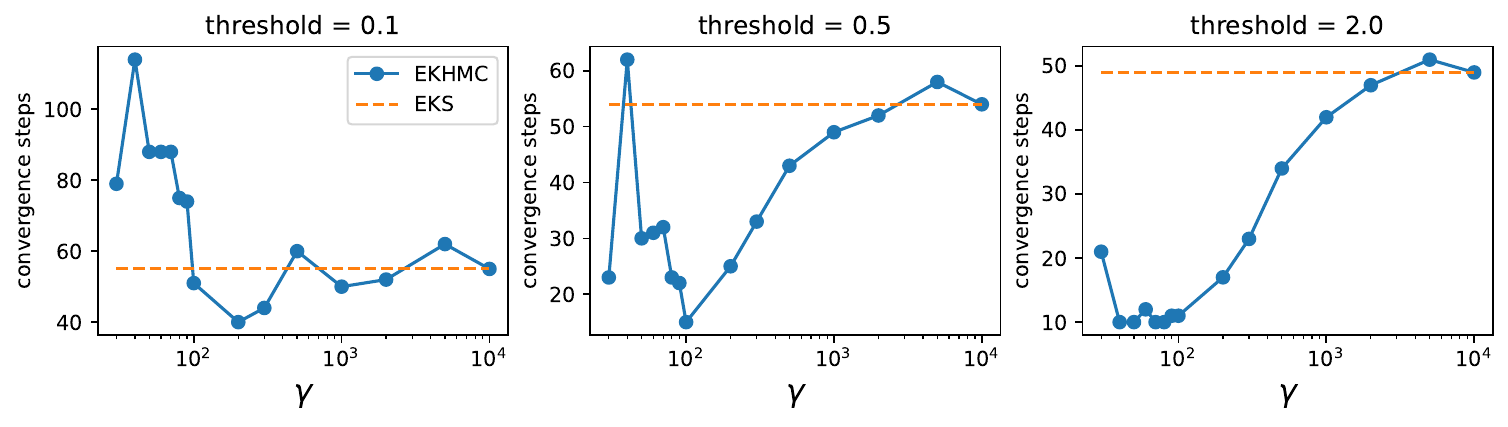}
    \caption{\yw{Convergence time (of $u_2$ mean) as a function of damping coefficient $\gamma$. Left, Middle, Right: thresholds are 0.1, 0.5, 2.0, respectively. The takeaway from these plots is: large damping converges faster eventually (small threshold), while small damping converges faster initially (large threshold).}}
    \label{fig:gamma}
\end{figure}

 We choose $a=0.01$, $\epsilon=0.2$ in both EKS and EKHMC. {We find 
choosing $\gamma=1$ causes overshooting \yw{(the trajectory exhibits 
oscillatory rather than monotonic convergence)}, a phenomenon which
can be ameliorated by increasing to \yw{$\gamma=100$}; indeed this latter
value appears, empirically, to be close to optimal in terms of convergence speed. We conjecture this difference from the linear case, where the
optimal value is $\gamma=1.83$ (see Remark \ref{r:5.2}),
is due to the non-Gaussianity of the desired target: the particles have accumulated considerable momentum when entering the linear convergence regime, and so extra damping is required to counteract this and avoid overshooting. Despite the
desirable problem independence of the optimal $\gamma$ in the linear case,
we expect case-by-case optimization to be needed for nonlinear
problems.} \yw{We also empirically study how $\gamma$ affects the convergence speed for this 2D problem: we sweep $\gamma\in [30,10^4]$, compute the number of steps needed to be close to the convergent point by a threshold, shown in FIG.~\ref{fig:gamma}. When the threshold is small, large $\gamma$ converges faster; when the threshold is large, small $\gamma$ converges faster. The takeaway is that: large $\gamma$ converges faster eventually (small threshold), while small $\gamma$ converges faster initially (large threshold).}

We evolve the ensemble of particles for 200 iterations, and record their positions in the last iteration as the approximation of the posterior distribution, {shown in Figure \ref{fig:2d}.} The EKHMC samples are quite similar to EKS samples, both of which are reasonable approximations of the samples obtained from the gold standard HMC \cite{duane1987hybrid} simulations -- but both approximations of the gold standard miss the full spread of the
true distribution, because of the ensemble Kalman approximation they invoke.
We also compute four {ensemble} quantities to characterize the evolution of EKHMC and EKS: $u_1$ mean, $u_2$ mean, two eigenvalues $\sigma_1$ and $\sigma_2$ of the covariance matrix $C_q(Z)$. EKHMC has faster convergence towards equilibrium than EKS, benefiting from the second-order dynamics.

\subsection{Darcy Flow}

\begin{figure}
    \centering
    \includegraphics[width=1.0\linewidth]{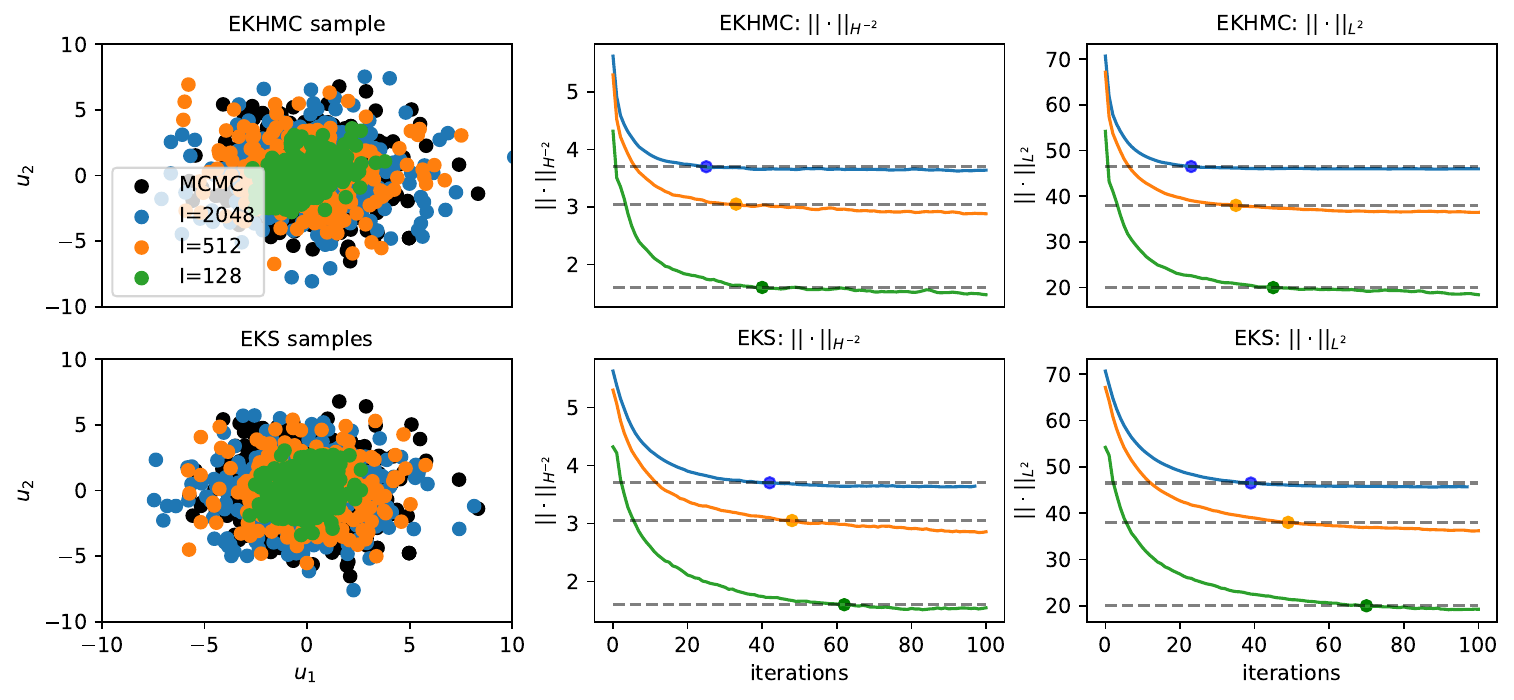}
    \caption{\yw{Darcy flow. Left: samples obtained from EKHMC (top) and EKS (bottom), compared with MCMC. Middle: Evolution of $||u||_{H^{-2}}$ for EKHMC and EKS for different $I=128,512,2048$. Right: The same as the middle, but $||u||_{L^2}$ instead of $||u||_{H^{-2}}$. EKHMC converges faster than EKS.}}
    \label{fig:darcy}
\end{figure}

This example follows the problem setting detailed in~\cite{garbuno2020interacting}. We summarize the essential problem specification here for completeness.{The forward
problem of porous medium flow, defined by permeability $a(\cdot)$ and source term $f(\cdot)$, is to find the pressure field $p(\cdot)$
where $p(\cdot)$ is solution of the following elliptic PDE:}
\begin{equation}
\begin{aligned}
    -\nabla\cdot \Bigl(a(x)\nabla p(x)\Bigr)&=f(x)\,, \quad x\in D=[0,1]^2\,, \\
    p(x)&=0\,, \quad x\in \partial D\,.
\end{aligned}
\end{equation}
We assume that the permeability $a(x)=a(x;u)$ depends on some unknown parameters $u\in\mathbb{R}^d$. 
{The resulting  inverse problem is, given (noisy) pointwise measurements of $p(x)$
on a grid, to infer $a(x;u)$ or $u$.}
We model $a(x;u)$ as a log-Gaussian field. The Gaussian underlying 
this log-Gaussian model has mean zero and has
precision operator defined as $\mathcal{C}^{-1}=(-\triangle+\tau^2\mathcal{I})^\alpha$; here
$\triangle$ is equipped with Neumann boundary conditions on the spatial-mean zero functions. We set $\tau=3$ and $\alpha=2$ in the experiments. Such parametrization yields as Karhunen-Loeve (KL) expansion

\begin{equation}
    {\rm log} a(x;u)=\sum_{l\in K} u_l\sqrt{\lambda_l}\varphi_l(x)\,,
\end{equation}
where $\varphi_l(x)={\rm cos}(\pi\langle l,x\rangle)$, $\lambda_l=(\pi^2|\ell|^2+\tau^2)^{-\alpha}$, $K\equiv \mathbb{Z}^2$. {Draws from this
random field are H\"older with any exponent less than one.}
In practice, we truncate $K$ to have dimension $d=2^8$. We generate a truth random field by sampling $u\sim \mathsf{N}(0,I_d)$. We create data $y$ with $\eta\sim \mathsf{N}(0,0.1^2 I_K)$. We choose the prior covariance $\Gamma_0=10^2 I_d$. 

To characterize the evolution of EKHMC and compare it with EKS, we compute two metrics:
\begin{equation}
    d_{H^{-2}}(\cdot)=\sqrt{\frac{1}{I}\sum_{j=1}^I  ||u^{(j)}(t)-\cdot||_{H^{-2}}^2}\,,\quad d_{L^2}(\cdot)=\sqrt{\frac{1}{I}\sum_{j=1}^I  ||u^{(j)}(t)-\cdot||^{2}_{L^2}}\,.
\end{equation}
For these metrics, we use norms 
\begin{equation}
    ||u||_{H^{-2}}=\sqrt{\sum_{l\in K_d}|u_l|^2\lambda_l}\,, \quad ||u||_{L^2}=\sqrt{\sum_{l\in K_d} |u_l|^2}\,,
\end{equation}
where the first is defined in the negative Sobolev space $H^{-2}$, whilst the second in the $L^2$ space.

We set $a=0.01$ and $\epsilon=1.0$ for both EKHMC and EKS. We set $\gamma=1$ in EKHMC; {unlike the previous example, this choice does not lead to over-shooting.}
We simulate the particles for $100$ iterations, and use their positions in the last iteration as the approximation to the posterior distribution. We compare the samples obtained from MCMC (we use the pCN variant on RWMH~\cite{cotter2013mcmc}), EKS ($I=128,512,2048$) and EKHMC ($I=128,512,2048$) in Figure~\ref{fig:darcy}. The evolution of $||u||_{H^{-2}}$ and $||u||_{L^2}$ are plotted to show that convergence is achieved, and that EKHMC converges faster than EKS.

\section{Conclusions}\label{5}

{Gradient-free methods for inverse problems are increasingly important in large-scale
multi-physics science and engineering problems; see \cite{huang2022efficient} and
the references therein. In this paper
we have provided an initial study of gradient-free methods which leverage the potential
power of Hamiltonian based sampling. Analysis of the resulting method is hard, and
our initial theoretical results leave many avenues open; in particular convergence to
equilibrium {is not established for the underlying nonlinear, nonlocal
Fokker-Planck equation arising from the mean field model,
and optimization of convergence rates over $\gamma$ is not understood,
except in the setting of linear
Gaussian inverse problems.}
The Fokker-Planck equation associated with standard underdamped Langevin dynamics has 
been studied in the context of hypocoercity -- see \cite{bakry1985diffusions,villani2006hypocoercivity} -- and generalization of these methods will be of potential interest. Preconditioned HMC is also proven to converge in \cite{bou2021two}; generalization to ensemble approximations would be of interest. On the computational side there are other approaches to avoiding gradient computations, yet leveraging Hamiltonian structure, which need to be evaluated and compared with what is proposed here \cite{lan2016emulation}.}

\section{Appendix}\label{6}
We present the proofs of {various results from the paper here.}

\subsection{Proof of Proposition \ref{prop.prop0}}
\label{proofprop0}
\begin{proof}
{The determination of the most general form of diffusion process which
is invariant with respect to a given measure has a long history; see
\cite{duncan2017using}[Theorem 1] for a statement 
and the historical context. Note, however, that this theory
is not developed in the mean-field setting and concerns only the
linear Fokker-Planck equation. However, the dependence of the matrices 
$\mathcal{K}, \mathcal{J}$ and $\mathcal{M}$ on $\rho$ readily allows use
of the approach taken in the linear case. We find it expedient to use
the exposition of this topic in \cite{ma2015complete}. The same
derivation as used to obtain Eq. (5) from \cite{ma2015complete} 
shows\footnote{We use the notational convention concerning divergence
of matrix fields that is standard in continuum mechanics \cite{gonzalez2008first};
this differs from the notational convention adopted in \cite{ma2015complete}.}
that the density $\rho$ associated with the dynamics \eqref{eqsol}, \eqref{eqsol2} satisfies Eq. \eqref{fp}; this follows from the (resp. skew-)
symmetry properties of (resp. $\mathcal{J}$) $\mathcal{K}$ and the fact
that $\mathcal{J}$, $\mathcal{K}$ and  $\mathcal{M}$ are assumed independent
of $z$, despite their dependence on $\rho.$ It is also manifestly the case
that the Gibbs measure is invariant for \eqref{fp} since the
mass term $\mathcal{M}$ appearing in $\mathcal{H}$ is independent of $z$ so that
$$\rho\nabla \cH+\nabla\rho=0$$
when $\rho$ is the Gibbs measure
and \eqref{fp} shows that then $\partial_t \rho=0.$}
\end{proof}

\subsection{Proof of Proposition \ref{aft}}
\label{proofaf}
\begin{proof}
In fact, consider the affine transformation 
\[q^{(i)}=A v^{(i)}+b\,,\quad p^{(i)}=A u^{(i)}\,,\] along with  $$\tilde{\Phi}(v^{(i)})=\Phi(q^{(i)})=\Phi(Av^{(i)}+b)\,,\quad\tilde{H}(v^{(i)}, u^{(i)})=H(q^{(i)}, p^{(i)})=H(Av^{(i)}+b,Au^{(i)})\,.$$ Then we have that the gradient term scales like $\nabla_{v^{(i)}}\tilde{\Phi}(v^{(i)})$= $A^T\nabla_{q^{(i)}}\Phi(q^{(i)})$ and the covariance preconditioner scales like $C_v={A^{-1}} C_q {A^{-1}}^T$. The generalized square root scales like $\sqrt{C_v}={A^{-1}} \sqrt{C_q}$. Therefore we can check that affine invariance holds for the particle systems \eqref{ensol} and \eqref{eqmom}. Affine invariance for the gradient free approximation  \eqref{encgf} is more easily seen to be true.  Finally for the Fokker-Planck equation \eqref{fp}, we can check similarly via the scaling of the terms. See a similar argument in the proof of Lemma 4.7 in \cite{garbuno2020affine}.
\end{proof}
\subsection{Proof of Proposition \ref{afli}}
\label{proof3.3}
\begin{proof}
{We can take expectations in \eqref{linearp} to obtain the first two ODEs
in \eqref{vare}. Now define $$\hat{z}=z-m(\rho)\triangleq\begin{pmatrix}
\hat{q} \\
\hat{p}
\end{pmatrix}\,,$$  to obtain the following evolution for $\hat{z}$:
\[\begin{aligned}
\dot{\hat{q}}&= \hat{p}\,, \\
\dot{\hat{p}}&=-\cC_q(\rho)B^{-1}\hat{q} - {\gamma}\hat{p}+\sqrt{2 {\gamma}\cC_q(\rho)}\dot{W}\,.
\end{aligned}\]
Note that $\cC_q=\mathop{\mathbb{E}}[\hat{q}\otimes \hat{q}]$, $\cC_p=\mathop{\mathbb{E}}[\hat{p}\otimes \hat{p}]$, and $\cC_{q,p}=\mathop{\mathbb{E}}[\hat{q}\otimes \hat{p}]$. We can 
use Ito's formula and the above evolution to derive the last three ODEs in \eqref{vare}. 
The form of the steady state solutions is immediate from setting the right hand side of
\eqref{vare} to zero.}

{Now, we will establish the independence of the essential dynamics of the
mean and covariance on the specific choice of $B,c$. To this end, denote $x_1 = B^{1/2}(B^{-1}m_q-c)$, $x_2 = B^{-1/2}m_p$, $X = B^{-1/2}\cC_qB^{-1/2}$, $Y = B^{-1/2}\cC_pB^{-1/2}$, $Z = B^{-1/2}\cC_{q,p}B^{-1/2}$. Applying this
change of variables we obtain, from \eqref{vare}, the system} 
\begin{equation}
\label{vares}
\begin{aligned}
\dot{x_1}&= x_2\,, \\
\dot{x_2}&=-X x_1 - {\gamma}x_2\,,\\
\dot{X}&=Z+Z^T\,,\\
\dot{Y}&=-X Z-Z^T X- 2{\gamma}Y+2{\gamma}X\,,\\
\dot{Z}&=- {\gamma}Z -X^2 +Y\,.
\end{aligned}
\end{equation}
We notice that the steady state solutions take the form $(x_1,x_2,X,Y,Z)=(w,0,X,X,0)$ for a projection matrix $X^2=X$ and  $w$ such that $Xw=0$.
{The  steady state} with $X=I_N$ and $w=0$ corresponds to the desired posterior mean.
The transformed equation is indeed independent of $B$ and $c$. Therefore the speed of convergence  $(x_1,x_2,X,Y,Z) \to (0,0,I_N,I_N,0)$, within its basin of
attraction, is indeed independent of $B$ and $c.$ The same is true in the original
variables.

{To complete the proof of stability it suffices} to establish the local exponential 
convergence to the
steady solution $(x_1,x_2,X,Y,Z)=(0,0,I_N,I_N,0)$. As a first step, we show that the following system converges to $(X,Y,Z)=(I_N,I_N,0)$:
\begin{equation}
\label{meanco}
\begin{aligned}
\dot{X}&=Z+Z^T\,,\\
\dot{Y}&=-X Z-Z^T X- 2{\gamma}Y+2{\gamma}X\,,\\
\dot{Z}&=-{\gamma}Z -X^2 +Y\,.
\end{aligned}
\end{equation}
{Linearization around $(X,Y,Z)=(I_N,I_N,0)$ {in variables of entries of $(X,Y,Z+Z^T)$} gives the matrix $$\begin{pmatrix}
0 & 0 & 2I_{N^2} \\
2\gamma I_{N^2} & -2\gamma I_{N^2} & -2I_{N^2} \\
-2I_{N^2} & I_{N^2} & -\gamma I_{N^2} 
\end{pmatrix}$$
{whose eigenvalues satisfy $x^3+3\gamma x^2+(2\gamma ^2+6)x+4\gamma=0$ which all have
strictly negative real part for $\gamma>0$, since we can show the unique real eigenvalue is in the interval $(-\gamma,0)$.} Therefore we conclude the local exponential convergence of covariances $X,Y,Z$ in a neighbourhood of $(X,Y,Z)=(I_N,I_N,0)$.
The exponential convergence of means $(x_1,x_2) \to (0,0)$ then follows 
linearizing the first two linear ODEs in the system \eqref{vares} at $X=I_N.$}

{Similarly, for any other steady state $(X,Y,Z)=(X,X,0)$ with a projection matrix $X\neq I_N$, we will show that there is an unstable direction in \eqref{meanco}. In fact, since $X$ is symmetric with only eigenvalues $1$ and $0$, there exists a nonzero vector $v$ such that $Xv=0$. Linearization around $(X,Y,Z)=(X,X,0)$ in the direction $(avv^T,bvv^T,cvv^T)$ gives the following $3\times3$ matrix for scalars $a,b,c$:
$$\begin{pmatrix}
0 & 0 & 2 \\
2\gamma  & -2\gamma  & 0 \\
0 & 1 & -\gamma
\end{pmatrix}$$
{whose eigenvalues satisfy $x^3+3\gamma x^2+2\gamma ^2x-4\gamma=0$. {For
all $\gamma>0$ it may be shown that there exists a real eigenvalue in the 
interval $(0,1)$, \yw{since $f(0)f(1)<0$}; thus there is an unstable direction determined by $(a,b,c)$,} and therefore in the original formulation.}} 
\end{proof}
\begin{remark}
\label{optg}
{We can further investigate the spectral gap around $(X,Y,Z)=(I_N,I_N,0)$, which is the absolute value in the real root of equation $x^3+3\gamma x^2+(2\gamma ^2+6)x+4\gamma=0$.
To be precise, we can show that for $x_0=-\sqrt{12-\sqrt{128}}$ and $\gamma_0=-(4+3x_0^2)/(4x_0)\approx 1.83$, the spectral gap is maximized as $-x_0$ and we expect fastest convergence for our method in the linear setting.}

{To establish this claim we proceed as follows. By the intermediate value theorem, in order to show there always exists a root in $[x_0,0)$ and the spectral gap is at most $-x_0$, we only need to show that $$x_0^3+3\gamma x_0^2+(2\gamma ^2+6)x_0+4\gamma\leq0\,.$$  The claim is true because $x_0<0$ and $2x_0(x_0^3+6x_0)=(4+3x_0^2)^2/4$. By the basic inequality $a^2+b^2\geq 2ab$, we have  $$x_0^3+3\gamma x_0^2+(2\gamma ^2+6)x_0+4\gamma=2x_0\gamma^2+(4+3x_0^2)\gamma+x_0^3+6x_0\leq0\,.$$
The maximal spectral gap is attained if and only if the equality in $a^2+b^2\geq 2ab$ holds, i.e. when $\gamma=\gamma_0$.}

We also plot the spectral gap as a function of $\gamma$ to help visualize the dependence of the rate of convergence on damping for the linear problem; see Figure \ref{fig:damp}. {The clear message from this figure is that there is a natural
optimization problem for $\gamma$ in this linear Gaussian setting; this is used
to motivate searches for optimal parameters in the non-Gaussian setting.}

\begin{figure}
    \centering
    \includegraphics[width=0.6\linewidth]{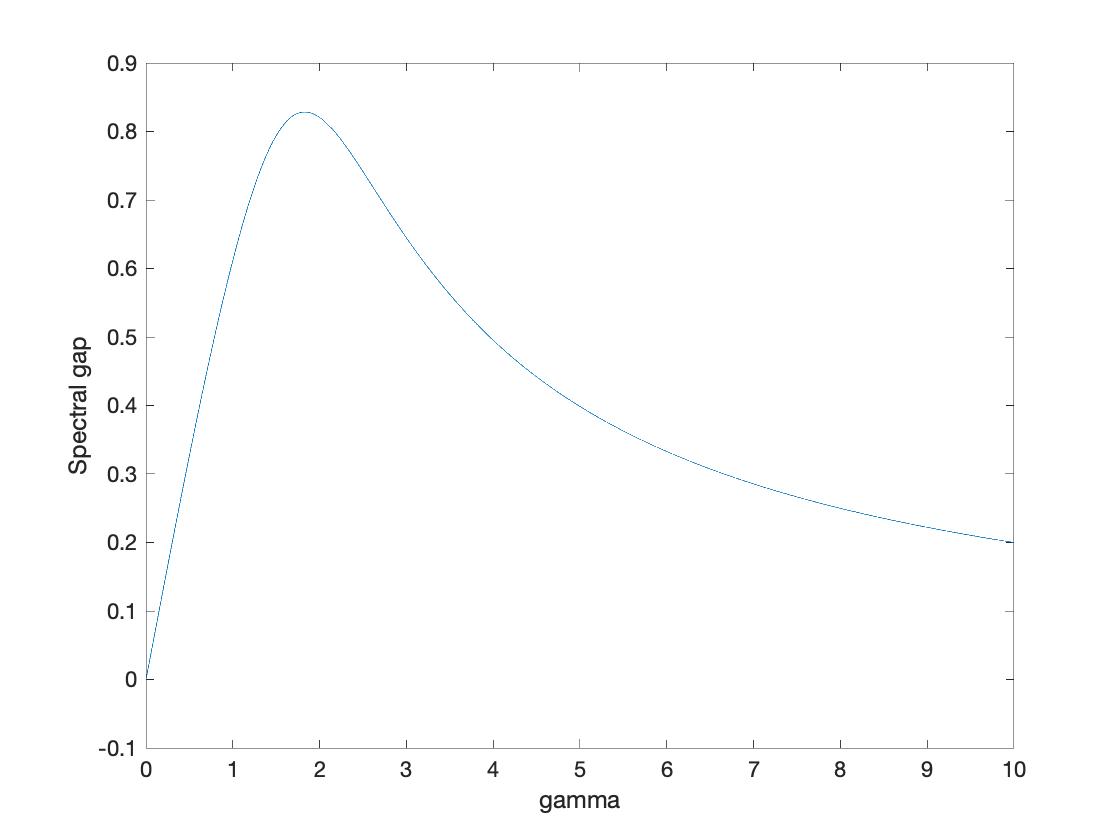}
    \caption{Spectral gap as a function of $\gamma$}
    \label{fig:damp}
\end{figure}
\end{remark}

\subsection{Proof of Proposition \ref{prop3}}
\label{proof3}
\begin{proof}
Note that since $m(\rho)$, $\cC(\rho)$ are independent of the particle $z$, we have
\[\nabla \rho=-{\cC}(\rho)^{-1}(z-m(\rho)) \rho\,,\]
\[D^2 \rho=\left(-{\cC}(\rho)^{-1}+\left({\cC}(\rho)^{-1}(z-m(\rho))\right) \otimes\left({\cC}(\rho)^{-1}(z-m(\rho))\right)\right) \rho\,.\]
{Using this we can substitute the Gaussian profile into equation \eqref{linearfp}. We compute the first term on the right hand to obtain}
\[\left(-{\cC}(\rho)^{-1}(z-m(\rho)) \cdot \begin{pmatrix}
-p \\
\cC_q(\rho)(B^{-1}q-c) + {\gamma}p
\end{pmatrix}+{\gamma}N\right) \rho \triangleq I_1\rho\,.\]
The second term on the right hand side can be written as:
\[\begin{pmatrix}
0 & -\cC_q(\rho) \\
\cC_q(\rho) & {\gamma}\cC_q(\rho)
\end{pmatrix}:\left(-{\cC}(\rho)^{-1}+\left({\cC}(\rho)^{-1}(z-m(\rho))\right) \otimes\left({\cC}(\rho)^{-1}(z-m(\rho))\right)\right) \rho\triangleq (I_2+I_3)\rho \,.\]
For the left hand side of \eqref{linearfp}, note that \[\frac{\mathrm{d}}{\mathrm{d} t} \operatorname{det} {\cC}(\rho)=\operatorname{Tr}\left[\operatorname{det} {\cC}(\rho) {\cC}(\rho)^{-1} \frac{\mathrm{d}}{\mathrm{d} t} {\cC}(\rho)\right], \,\, \frac{\mathrm{d}}{\mathrm{d} t} {\cC}(\rho)^{-1}=-{\cC}(\rho)^{-1}\left(\frac{\mathrm{d}}{\mathrm{d} t} {\cC}(\rho)\right) {\cC}(\rho)^{-1}\,.\]
{Using the evolution of the mean \eqref{vare} we obtain}
\[\begin{aligned}
\frac{\mathrm{d}}{\mathrm{d} t}\|z-m(\rho)\|_{{\cC}(\rho)}^{2}=& 2\left\langle\frac{\mathrm{d}}{\mathrm{d} t}(z-m(\rho)), {\cC}(\rho)^{-1}(z-m(\rho))\right\rangle \\
&+\left\langle(z-m(\rho)), \frac{\mathrm{d}}{\mathrm{d} t}\left({\cC}(\rho)^{-1}\right)(z-m(\rho))\right\rangle \\
=&2\left\langle \begin{pmatrix}
-m_p(\rho) \\
\cC_q(\rho)(B^{-1}m_q(\rho)-c) + {\gamma}m_p(\rho)
\end{pmatrix}, \cC(\rho)^{-1}(z-m(\rho))\right\rangle\\&-\left\langle{\cC}(\rho)^{-1}(z-m(\rho)), \frac{\mathrm{d}}{\mathrm{d} t}\left({\cC}(\rho)\right){\cC}(\rho)^{-1}(z-m(\rho))\right\rangle\triangleq 2(I_4-I_5)\,.
\end{aligned}\]
Therefore we can compute the left hand side as:\[\begin{aligned}
\partial_{t} \rho &=\left[-\frac{1}{2}(\operatorname{det} {\cC}(\rho))^{-1}\left(\frac{\mathrm{d}}{\mathrm{d} t} \operatorname{det} {\cC}(\rho)\right)-\frac{1}{2} \frac{\mathrm{d}}{\mathrm{d} t}\left\|z-m(\rho)\right\|_{{\cC}(\rho)}^{2}\right] \rho \\
&=\left[-\frac{1}{2}\operatorname{Tr}\left[{\cC}(\rho)^{-1} \frac{\mathrm{d}}{\mathrm{d} t} {\cC}(\rho)\right]-I_4+I_5\right] \rho\triangleq (I_6-I_4+I_5)\rho\,.
\end{aligned}\]

In order to show that the Gaussian profile is indeed a solution to \eqref{linearfp}, we only need to show that $I_1+I_2+I_3=-I_4+I_5+I_6$. Note that \[\begin{aligned}I_1+I_4&={\gamma}N+{\cC}(\rho)^{-1}(z-m(\rho)) \cdot \begin{pmatrix}
p-m_p(\rho) \\
-\cC_q(\rho)B^{-1}(q-m_q(\rho)) - {\gamma}(p-m_p(\rho))
\end{pmatrix}\\&={\gamma}N+{\cC}(\rho)^{-1}(z-m(\rho)) \cdot\begin{pmatrix}
0 &I_N \\
-\cC_q(\rho)B^{-1} &- {\gamma}I_N
\end{pmatrix} (z-m(\rho))\,.\end{aligned}\]
Defining $\hat{z}=z-m(\rho)$, we collect the terms in the to-be-proven identity $I_1+I_2+I_3=-I_4+I_5+I_6$ as:
\[\begin{aligned}
&\hat{z}^T\left[{\cC}(\rho)^{-1}\begin{pmatrix}
0 &I_N \\
-\cC_q(\rho)B^{-1} &- {\gamma}I_N
\end{pmatrix}+{\cC}(\rho)^{-1} \left(\begin{pmatrix}
0 & -\cC_q(\rho) \\
\cC_q(\rho) & {\gamma}\cC_q(\rho)
\end{pmatrix}-\frac{1}{2}\frac{\mathrm{d}}{\mathrm{d} t}\left({\cC}(\rho)\right)\right){\cC}(\rho)^{-1}\right]\hat{z}\\+&{\gamma}N-\operatorname{Tr}\left[\begin{pmatrix}
0 & -\cC_q(\rho) \\
\cC_q(\rho) & {\gamma}\cC_q(\rho)
\end{pmatrix}{\cC}(\rho)^{-1}\right]+\frac{1}{2}\operatorname{Tr}\left[\frac{\mathrm{d}}{\mathrm{d} t} {\cC}(\rho){\cC}(\rho)^{-1} \right]\triangleq \hat{z}^T M_1\hat{z}+M_2=0\,.
\end{aligned}\]
We only need to show that $M_1+M_1^T=M_2=0$.

In fact, we compute ${\cC}(\rho)(M_1+M_1^T){\cC}(\rho)$ to be:
\[\begin{pmatrix}
0 &I_N \\
-\cC_q(\rho)B^{-1} &- {\gamma}I_N
\end{pmatrix}{\cC}(\rho)+{\cC}(\rho)\begin{pmatrix}
0 &-(B^{-1})^T\cC_q(\rho) \\
I_N&- {\gamma}I_N
\end{pmatrix}+\begin{pmatrix}
0 & 0 \\
0 & 2{\gamma}\cC_q(\rho)
\end{pmatrix}-\frac{\mathrm{d}}{\mathrm{d} t}\left({\cC}(\rho)\right)\,,\]
which equals zero upon blockwise computation via \eqref{vare}. Thus $M_1+M_1^T=0$.
Moreover, note that $M_2=-\operatorname{Tr}\left[\cC(\rho)M_1\right]$: $M_1+M_1^T=0$ implies $\cC(\rho)(M_1+M_1^T)=0$ and taking the trace on both sides {leads to the
conclusion} that $M_2=0$.
Therefore the Gaussian profile indeed solves the Fokker-Planck equation \eqref{linearfp}.
\end{proof}

\noindent{\bf Acknowledgments} The work of ZL is supported by IAIFI through NSF grant PHY-2019786. The work of AMS is supported by NSF award AGS1835860, the Office of 
Naval Research award N00014-17-1-2079 and by a Department of Defense Vannevar Bush Faculty Fellowship.

\bibliographystyle{abbrv}
\bibliography{ref_En}
          \end{document}